\documentclass[a4paper,11pt,twoside,english]{article}
\usepackage{babel}
\usepackage{latexsym,amsfonts}
\usepackage{amsthm}
\usepackage{amsmath}
\usepackage{mathtools}
\usepackage{paralist}
\usepackage{ifthen}

\usepackage{lastpage}
\usepackage{fancyhdr}
\fancypagestyle{plain}{

  \fancyhead[L,C]{}
  \fancyhead[R]{\today}
  \fancyfoot[R,L]{}
  \fancyfoot[C]{\thepage \ of \pageref{LastPage}}
  }

\title{A remark on condensation of singularities}
\author{Jan-David Hardtke}
\date{}

\providecommand{\ssq}{\subseteq}

\providecommand{\N}{\ensuremath{\mathbb{N}}}

\providecommand{\R}{\ensuremath{\mathbb{R}}}
\providecommand{\C}{\ensuremath{\mathbb{C}}}
\providecommand{\K}{\ensuremath{\mathbb{K}}}

\providecommand{\eps}{\ensuremath{\varepsilon}}
\providecommand{\Norm}[1]{\ensuremath{\left|\hspace{-1.2pt}\left|\hspace{-1.2pt}\left|#1\right|\hspace{-1.2pt}\right|\hspace{-1.2pt}\right|}}

\providecommand{\keywords}[1]{
  {\let\thefootnote=\relax
  \footnote{{\em Keywords}: #1}}
  \addtocounter{footnote}{-1}
  }

\providecommand{\AMS}[1]{
  {\let\thefootnote=\relax
  \footnote{{\em AMS Subject Classification} (2010): #1}}
  \addtocounter{footnote}{-1}
  }

\providecommand{\address}{
  {\sc \noindent Department of Mathematics \\
  Freie Universit\"at Berlin \\
  Arnimallee 6, 14195 berlin \\
  Germany \\}
  }

\DeclarePairedDelimiter{\set}{\lbrace}{\rbrace}
\DeclarePairedDelimiter{\paren}{\lparen}{\rparen}

\DeclarePairedDelimiter{\abs}{\lvert}{\rvert}
\DeclarePairedDelimiter{\norm}{\lVert}{\rVert}

\theoremstyle{definition}

\theoremstyle{plain}

\newtheorem*{lemma*}{Lemma}

\newtheorem*{theorem*}{Theorem}

\newenvironment{Proof}[1][\proofname]{\begin{proof}[#1] \setlength{\parindent}{0pt}}{\end{proof}}
\newenvironment{Abstract}{\centering\begin{minipage}{0.8\textwidth} \noindent \small {\sc Abstract.}}{\end{minipage}\par}

\usepackage{color}
\definecolor{darkgreen}{rgb}{0,0.5,0}

\newtagform{colored}[\color{blue}]{\color{blue}(}{\color{blue})}
\usetagform{colored}

\usepackage[colorlinks,linkcolor=blue,citecolor=red,urlcolor=darkgreen]{hyperref}
\providecommand{\email}{{\it E-mail address:} \href{mailto:hardtke@math.fu-berlin.de}{\tt hardtke@math.fu-berlin.de}}

\usepackage{amsrefs}

\begin{document}

\maketitle

\begin{Abstract}
Recently in \cite{sokal}  Alan D. Sokal gave a very short and completely elementary proof of 
the uniform boundedness principle. The aim of this note is to point out that by using a similiar
technique one can give a considerably short and simple proof of a stronger statement, namely a  
principle of condensation of singularities for certain double-sequences of non-linear operators  
on quasi-Banach spaces, which is a bit more general than a result of I.\,S. G\'al from \cite{gal}.
\end{Abstract}
\AMS{46A16; 47H99}

\ \par
Let us begin by recalling that a quasi-norm on a linear space $X$ over the field $\K=\R$ or $\C$ is 
a map $\norm{\,.\,}:X \to [0,\infty)$ such that
\begin{enumerate}[(i)]
\item $\norm{x}=0 \ \Rightarrow \ x=0,$
\item $\norm{\lambda x}=\abs{\lambda}\norm{x} \ \ \forall x\in X, \forall \lambda\in \K,$
\item $\exists K\geq1$ such that $\norm{x+y}\leq K(\norm{x}+\norm{y}) \ \ \forall x,y\in X$.
\end{enumerate}
The least constant $K\geq1$ which fulfils (iii) is sometimes called the modulus of concavity of 
the quasi-norm and the pair $(X,\norm{\,.\,})$ is called a quasi-normed space.\par
If $p\in (0,1]$ then a map $\norm{\,.\,}:X \to [0,\infty)$ is called a $p$-norm if it satisfies (i) and (ii)
and the condition
\begin{enumerate}[(iv)]
\item $\norm{x+y}^p\leq \norm{x}^p+\norm{y}^p \ \ \forall x,y\in X$.
\end{enumerate}
The pair $(X,\norm{\,.\,})$ is then called a $p$-normed space. Every $p$-norm is a quasi-norm with $K=2^{1/p-1}$. 
The standard examples of $p$-normed spaces are of course the spaces $L^p(\mu)$ with the $p$-norm 
\begin{equation*}
\norm{f}_p=\paren*{\int_{S}\,\abs{f}^p\,\mathrm{d}\mu}^{1/p}
\end{equation*}
for any measure space $(S,\mathcal{A},\mu)$, which includes in particular the sequence spaces $\ell^p$.\par
The Aoki-Rolewicz-Theorem (cf. \cites{aoki,rolewicz} or \cite{kalton2}*{Lemma 1.1}) states that if 
$\norm{\,.\,}$ is a quasi-norm on $X$ with modulus of concavity $\leq K$ and if $p\in (0,1]$ is defined 
by $(2K)^p=2$ then for any $x_1,\dots,x_n\in X$ the inequality
\begin{equation*}
\norm*{\sum_{i=1}^nx_i}^p\leq 4\sum_{i=1}^n\norm{x_i}^p
\end{equation*}
holds.\par
Then an equivalent $p$-norm on $X$ can be defined by
\begin{equation*}
\Norm{x}=\inf\set*{\paren*{\sum_{i=1}^n\norm{x_i}^p}^{1/p}:n\in \N,\,x_1,\dots,x_n\in X,\,x=\sum_{i=1}^nx_i},
\end{equation*}
cf. \cites{aoki,rolewicz} or \cite{kalton2}*{Theorems 1.2 and 1.3}.\par
The topology $\tau$ on $X$ induced by the quasi-norm $\norm{\,.\,}$ is defined by declaring a subset $O\ssq X$ 
to be open if for every $x\in O$ there exists $\eps>0$ such that $\set*{y\in X: \norm{x-y}<\eps}\ssq O$, which 
is the same as the topology induced by the metric $d$ that is defined by $d(x,y)=\Norm{x-y}^p$. The space $(X,\tau)$
is a topological vectorspace and $(X,\norm{\,.\,})$ is called a quasi-Banach space if $(X,\tau)$ is complete (or
equivalently if $X$ is complete with respect to the above metric $d$ or any other translation-invariant metric on 
$X$ that induces the topology $\tau$).\par
A linear operator $T$ between two quasi-normed spaces $X$ and $Y$ is continuous if and only if there is a constant $M\geq 0$
such that $\norm{Tx}\leq M\norm{x}$ for every $x\in X$. The space of continuous linear operators from $X$ into $Y$
is denoted by $L(X,Y)$. A quasi-norm on $L(X,Y)$ can be defined by $\norm{T}=\sup\set{\norm{Tx}:x\in X,\,\norm{x}\leq1}$.
The modulus of concavity of this quasi-norm is less than or equal to the one of $Y$ and $L(X,Y)$ is complete if $Y$ is complete, 
in particular the dual space $X^*=L(X,\K)$ is always a Banach space.\par
It should be mentioned that the Hahn-Banach-Theorem fails in general for quasi-normed spaces, some of them do not even
possess a separating dual, for example it is a well-known result that for $p\in (0,1)$ the dual space of $L^p[0,1]$ 
contains only $0$. However, the classical results on Banach spaces like the open mapping and the closed graph theorem, 
as well as the uniform boundedness principle still hold for quasi-Banach spaces with analogous proofs using the Baire
category theorem.\footnote{For more information on quasi-Banach spaces the reader is referred to \cite{kalton1} and references 
therein.}\par
The original proofs of the uniform boundedness principle by Banach and Hahn (cf. \cite{banach} and \cite{hahn}) made no use 
of the Baire category theorem, but of a technique which has come to be known as `gliding hump' method. The Baire category argument 
was originally introduced into this theory by S. Saks.\par
In 1927 Banach and Steinhaus proved a by now well-known generalisation of the uniform boundedness principle, the so called 
principle of condensation of singularities, which reads as follows.
\begin{theorem*}[Principle of condensation of singularities, cf. \cite{banachsteinhaus}]
Let $X$ be a Banach space and $(Y_n)_{n\in \N}$ a sequence of normed linear spaces. If $G_n\ssq L(X,Y_n)$
is unbounded for every $n\in \N$ then there exists an element $x\in X$ such that
\begin{equation*}
\sup_{T\in G_n}\norm{Tx}=\infty \ \ \forall n\in \N.
\end{equation*}
\end{theorem*}
\noindent The proof from \cite{banachsteinhaus} also uses the Baire theorem.\par
Over the years many generalisations of the unifrom boundedness principle and the principle of condensation of
singularities to suitable classes of non-linear operators on Banach and quasi-Banach spaces have been investigated.
In \cite{gal0} I.\,S. G\'al proved a principle of uniform boundedness for sequences of bounded homogenous operators 
(see below for definitions) on Banach spaces that satisfy certain conditions of `asymptotic subadditivity' and later in \cite{gal} 
he generalized his results to a principle of condensation of singularities for double-sequences of the same type. The authors of 
\cite{dickmeis}  consider quantitative versions of the principle of condensation of singularities for double-sequences of operators 
on quasi-Banach spaces that fulfil similiar conditions.\par
The Baire category theorem seems not to be useful to prove these non-linear (or even quantitave) extensions of the classical 
theorems, so the proofs use rather complicated variations and refinements of the original `gliding hump' method.\par
In \cite{sokal} Alan D. Sokal gave a very short and elegant proof of the classical uniform boundedness principle for linear 
operators on Banach spaces, which is also completely elementary (in particular, it does not use the Baire category theorem).\par
We wish to point out here that by using a suitable refinement of Sokal's idea we can give a proof of a theorem slightly
more general than the aforementioned result of G\'al that is considerably shorter than the original proof from \cite{gal}.
The technique we will use in the proof also resembles the one from \cite{dickmeis} but since we are not interested in 
quantitative versions it is less complicated.\par
To formulate the result we first recall that a map $T:X \rightarrow Y$, where $X$ and $Y$ are of course quasi-normed spaces,
is called homogenous if 
\begin{equation*}
\norm{T(\lambda x)}=\abs{\lambda}\norm{Tx} \ \ \forall x\in X, \forall \lambda \in \K
\end{equation*}
and it is called bounded if $\sup\set{\norm{Tx}:x\in X,\,\norm{x}\leq1}$ is finite. For a bounded homogenous map we denote
this supremum again by $\norm{T}$. It follows that $\norm{Tx}\leq \norm{T}\norm{x}$ for every $x\in X$ and $\norm{T}$ is the 
best possible such constant.\par
The result then reads as follows.
\begin{theorem*}
Let $X$ be a quasi-Banach space and $(Y_{nm})_{(n,m)\in \N^2}$ a double-sequence of quasi-normed spaces, as well as 
$(T_{nm}:X \rightarrow Y_{nm})_{(n,m)\in \N^2}$ a double-sequence of bounded homogenous operators satisfying
\begin{equation*}\label{+}
\sup_{n\in \N}\,\norm{T_{nm}}=\infty \ \ \forall m\in \N. \tag{+}
\end{equation*}
Suppose further that there are two sequences $(C_m)_{m\in \N}$ and $(L_m)_{m\in \N}$ of positive real numbers, a sequence 
$(N_m)_{m\in \N}$ in $\N$ and two double-sequences $(c_{nm})_{(n,m)\in \N^2}$ and $(f_{nm})_{(n,m)\in \N^2}$ of functions 
from $X$ into $[0,\infty)$ such that
\begin{enumerate}[\upshape(i)]
\item \label{i} for every $x\in X$, each $y\in X$ with $\norm{y}\leq 1$ and all $n,m\in \N$ with $n\geq N_m$ we have
\begin{equation*}
\norm{T_{nm}(x+y)}\leq C_m\paren*{\norm{T_{nm}x}+\norm{T_{nm}}\norm{y}+f_{nm}(x)},
\end{equation*}
\item \label{ii}for every $x\in X$ we have $f_{nm}(x)=O(1)$ as $n\to \infty$ uniformly in $m\in \N$,
\item \label{iii}for every $x\in X$, each $y\in X$ with $\norm{y}\leq 1$ and all $n,m\in \N$ the inequality
\begin{equation*}
\norm{T_{nm}y}\leq L_m\paren*{\norm{T_{nm}(x+y)}+\norm{T_{nm}x}+c_{nm}(x)\norm{T_{nm}}}
\end{equation*}
holds,
\item \label{iv}for every $x\in X$ and every $m\in \N$ we have $c_{nm}(x)\to 0$ as $n\to \infty$.
\end{enumerate}
Then there is an element $x\in X$ such that
\begin{equation*}
\sup_{n\in \N}\,\norm{T_{nm}x}=\infty \ \ \forall m\in \N.
\end{equation*}
\end{theorem*}
As we said, the conditions (i)--(iv) include in particular the case of G\'al's asymptotically subadditive double-sequences
(essentially only the functions $f_{nm}$ are new here). For some examples were such conditions of asymptotic subadditivity
occur naturally (e.\,g. the so called metric-mean interpolations) we refer the reader to \cite{gal0} and \cite{gal}.\par
Before we can come to the main proof we need an easy lemma.
\begin{lemma*}
For every sequence $(\alpha_n)_{n\in \N}$ of positive real numbers and every $\beta>0$ there is a sequence $(\beta_n)_{n\in \N}$
of positive real numbers such that
\begin{equation*}
\sum_{i=n+1}^{\infty}\beta_i<\alpha_n\beta_n \ \mathrm{and} \ \beta_n\leq \beta \ \ \forall n\in \N.
\end{equation*}
\end{lemma*}

\begin{Proof}
Just choose inductively $0<\beta_n\leq\beta$ such that  
\begin{equation*}
\sum_{i=m+1}^n\beta_i<\frac{\alpha_m\beta_m}{2} \ \ \forall m\in\set{1,\dots,n-1}.
\end{equation*}
\end{Proof}

\begin{Proof}[Proof of the theorem]
By the Aoki-Rolewicz-Theorem we may assume without loss of generality that $X$ is $p$-normed for some $p\in (0,1]$ (actually 
this is not necessary for the proof but it is more convenient).\par
Define $\psi: \N \rightarrow \N$ by $(\psi(1),\psi(2),\dots)=(1,1,2,1,2,3,1,2,3,4,\dots)$.\par
It easily follows from \eqref{iii} that for any $r\in (0,1]$, each $x\in X$ and all $n\geq N_m$ we have
\begin{align}\label{eq:1}
&\nonumber\frac{r\norm{T_{nm}}}{L_m}\leq \smashoperator[l]{\sup_{y\in B_r(x)}}\!\norm{T_{nm}y}+\norm{T_{nm}x}+c_{nm}(x)\norm{T_{nm}} \\
&\leq 2\smashoperator[l]{\sup_{y\in B_r(x)}}\!\norm{T_{nm}y}+c_{nm}(x)\norm{T_{nm}}, 
\end{align}
where $B_r(x)=\set*{y\in X: \norm{x-y}\leq r}$ (this is a kind of analogue of the lemma from \cite{sokal}).\par
According to the above lemma we can find a sequence $(\tilde{\beta}_n)_{n\in \N}$ in $(0,1]$ such that
\begin{equation}\label{eq:2}
\sum_{i=n+1}^{\infty}\tilde{\beta}_i<\frac{\tilde{\beta}_n}{8^pL_{\psi(n)}^pC_{\psi(n)}^p} \ \ \forall n\in \N.
\end{equation}
For every $n\in \N$ we put $\beta_n=\tilde{\beta}_n^{1/p}$ and 
\begin{equation}\label{eq:3}
\gamma_n=\frac{\beta_n}{8L_{\psi(n)}C_{\psi(n)}}-\paren*{\sum_{i=n+1}^{\infty}\beta_i^p}^{1/p},
\end{equation}
which by \eqref{eq:2} is strictly positive. Let us also put $x_0=0$.\par
Because of \eqref{iv} and our assumption \eqref{+} we can find $n_1\geq N_{\psi(1)}$ such that $c_{n_1\psi(1)}(x_0)\leq\beta_{1}/2L_{\psi(1)}$ 
and $\norm{T_{n_1\psi(1)}}\geq 1/\gamma_1$.\par
It follows from \eqref{eq:1} that 
\begin{equation*}
\frac{\beta_1}{4L_{\psi(1)}}\norm{T_{n_1\psi(1)}}\leq \smashoperator{\sup_{y\in B_{\beta_1}(x_0)}}\ \norm{T_{n_1\psi(1)}y}
\end{equation*}
and hence there is some $x_1\in X$ such that $\norm{x_1-x_0}\leq \beta_1$ and 
$\norm{T_{n_1\psi(1)}x_1}\geq \beta_1\norm{T_{n_1\psi(1)}}/8L_{\psi(1)}$.\par
Next we use \eqref{iv} and \eqref{+} to choose an index $n_2>\max\set{n_1,N_{\psi(2)}}$ such that $c_{n_2\psi(2)}(x_1)\leq \beta_2/2L_{\psi(2)}$ and 
$\norm{T_{n_2\psi(2)}}\geq 2/\gamma_2$ and then find (using \eqref{eq:1}) an element $x_2\in X$ with $\norm{x_2-x_1}\leq \beta_2$ and 
$\norm{T_{n_2\psi(2)}x_2}\geq \beta_2\norm{T_{n_2\psi(2)}}/8L_{\psi(2)}$.\par
Continuing in this way we obtain a strictly increasing sequence $(n_k)_{k\in \N}$ in $\N$ and a sequence $(x_k)_{k\in \N}$ in $X$ such that 
for every $k\in \N$ we have $n_k\geq N_{\psi(k)}$ and 
\begin{align}
&\norm{x_k-x_{k-1}}\leq \beta_k, \label{eq:4} \\
&\norm{T_{n_k\psi(k)}x_k}\geq \frac{\beta_k}{8L_{\psi(k)}}\norm{T_{n_k\psi(k)}}, \label{eq:5} \\
&\norm{T_{n_k\psi(k)}}\geq \frac{k}{\gamma_k}. \label{eq:6}
\end{align}
From \eqref{eq:4} it follows that $d(x_n,x_m)=\norm{x_n-x_m}^p\leq \sum_{i=m+1}^n\beta_i^p$ for every $n>m$ and hence
$(x_k)_{k\in \N}$ is a Cauchy sequence. By completeness, the limit $x=\lim_{k\to \infty}x_k$ exists and it follows that
\begin{equation}\label{eq:7}
\norm{x-x_m}^p\leq \sum_{i=m+1}^\infty\beta_i^p \ \ \forall m\in \N.
\end{equation}
Now for sufficiently large $k\in \N$ we have by \eqref{i}, \eqref{eq:5}, \eqref{eq:7}, \eqref{eq:3} and \eqref{eq:6}
\begin{align*}
&f_{n_k\psi(k)}(x)+\norm{T_{n_k\psi(k)}x}\geq \frac{\norm{T_{n_k\psi(k)}x_k}}{C_{\psi(k)}}-\norm{T_{n_k\psi(k)}}\norm{x_k-x} \\
&\geq \norm{T_{n_k\psi(k)}}\paren*{\frac{\beta_k}{8L_{\psi(k)}C_{\psi(k)}}-\paren*{\sum_{i=k+1}^{\infty}\beta_i^p}^{1/p}}=\norm{T_{n_k\psi(k)}}\gamma_k\geq k.
\end{align*}
Together with \eqref{ii} this implies $\norm{T_{n_k\psi(k)}x}\to \infty$ for $k\to \infty$. Since for every $m\in \N$ the set $\psi^{-1}(\set{m})$
is infinite it follows that $\sup_{n\in \N}\norm{T_{nm}x}=\infty$ for every $m\in \N$.
\end{Proof}

\address
\email

\end{document}